\documentclass[11pt]{article}

\usepackage{amscd,amsmath, amssymb, fancyhdr, epsfig, color, tocloft, mathtools}
\usepackage{dutchcal}

\usepackage[backref=page]{hyperref}
\renewcommand*{\backrefalt}[4]{%
	\ifcase #1 (Not cited.)%
	\or        (Cited on page~#2.)%
	\else      (Cited on pages~#2.)%
	\fi}

\hypersetup{
	colorlinks   = true,
	citecolor    = magenta,
	linkcolor    = blue,
	urlcolor     = magenta	
}

\newcommand{\version}{version 1.2,\ \ Feb. 5, 2023}
\makeatletter
\def\x@arrow{\DOTSB\Relbar}
\def\xlongequalsignfill@{\arrowfill@\x@arrow\Relbar\x@arrow}
\providecommand{\xlongequal}[2][]{%
	\ext@arrow 0099\xlongequalsignfill@{#1}{#2}}
\def\xlongrightarrowfill@{\arrowfill@\relbar\relbar\longrightarrow}
\makeatother

\numberwithin{equation}{section}

\def\eqref#1{(\ref{#1})}

\newcommand{\C}{{\mathbb C}}
\newcommand{\R}{{\mathbb R}}

\def\1{\sqrt{-1}\:}
\newcommand{\restrict}[1]{{\left|_{{\phantom{|}\!\!}_{#1}}\right.}}
\newcommand{\cntrct}                
{\hspace{2pt}\raisebox{1pt}{\text{$\lrcorner$}}\hspace{2pt}}
\newcommand{\arrow}{{\:\longrightarrow\:}}

\renewcommand{\phi}{\varphi}
\renewcommand{\epsilon}{\varepsilon}
\renewcommand{\geq}{\geqslant}

\newcommand{\Vol}{\operatorname{Vol}}

\newcommand{\GL}{\operatorname{GL}}

\newcounter{Mycounter}[section]
\newcounter{lemma}[section]
\newcounter{claim}[section]
\newcounter{sublemma}[section]
\newcounter{corollary}[section]
\newcounter{theorem}[section]
\renewcommand{\thetheorem}{{Theorem \thesection.\arabic{theorem}}}
\newcommand{\theorem}{%
	\setcounter{theorem}{\value{Mycounter}}
	\refstepcounter{theorem}
	\stepcounter{Mycounter}
	{\noindent \bf \thetheorem:\ }}

\newcounter{conjecture}[section]
\newcounter{proposition}[section]
\renewcommand{\theproposition} {{Proposition \thesection.\arabic{proposition}}}
\newcommand{\proposition}{%
	\setcounter{proposition}{\value{Mycounter}}
	\refstepcounter{proposition}
	\stepcounter{Mycounter}
	{\noindent \bf \theproposition:\ }}

\newcounter{definition}[section]
\renewcommand{\thedefinition} {{Definition~\thesection.\arabic{definition}}}
\newcommand{\definition}{%
	\setcounter{definition}{\value{Mycounter}}
	\refstepcounter{definition}
	\stepcounter{Mycounter}
	{\noindent \bf \thedefinition:\ }}

\newcounter{example}[section]
\renewcommand{\theexample}{{Example \thesection.\arabic{example}}}
\newcommand{\example}{%
	\setcounter{example}{\value{Mycounter}}
	\refstepcounter{example}
	\stepcounter{Mycounter}
	{\noindent \bf \theexample:\ }}

\newcounter{remark}[section]
\renewcommand{\theremark}{{Remark \thesection.\arabic{remark}}}
\newcommand{\remark}{%
	\setcounter{remark}{\value{Mycounter}}
	\refstepcounter{remark}
	\stepcounter{Mycounter}
	{\noindent \bf \theremark:\ }}

\newcounter{problem}[section]
\newcounter{question}[section]
\makeatletter

\@addtoreset{equation}{section}
\@addtoreset{footnote}{section}
\makeatother

\def\blacksquare{\hbox{\vrule width 5pt height 5pt depth 0pt}}
\def\endproof{\blacksquare}

\newcommand{\proof}{{\bf Proof: \ }}
\newcommand{\pstep}{{\bf Proof. Step 1: \ }}

\begin{document}
	
	\begin{center}
		{\Large\bf  Bimeromorphic geometry of 
LCK manifolds}\\[5mm]
		{\large
			Liviu Ornea\footnote{Liviu Ornea is  partially supported by Romanian Ministry of Education and Research, Program PN-III, Project number PN-III-P4-ID-PCE-2020-0025, Contract  30/04.02.2021},  
			Misha Verbitsky\footnote{Misha Verbitsky is partially supported 
				by the HSE University Basic Research Program,
FAPERJ 	SEI-260003/000410/2023 and CNPq - Process 310952/2021-2.\\[1mm]
				
			}\\[4mm]
			
		}
		
	\end{center}

	\hfill
	
	{\small
\hspace{0.15\linewidth}
\begin{minipage}[t]{0.7\linewidth}
{\bf Abstract} \\
A locally conformally K\"ahler (LCK) manifold
is a complex manifold $M$ which has a K\"ahler structure
on its cover, such that the deck transform group 
acts on it by homotheties. Assume that the K\"ahler
form is exact on the minimal K\"ahler cover of $M$. 
We prove that any bimeromorphic map $M'\rightarrow M$ 
is in fact holomorphic; in other words,
$M$ has a unique minimal model.
This can be applied to a wide class
of LCK manifolds, such as the Hopf manifolds,
their complex submanifolds and to OT manifolds.  
\end{minipage} 
	}
	\tableofcontents
	
\section{Introduction}

The notion of a minimal model is central 
in birational algebraic geometry since the XIX-th
centiry.  In 1901, Castelnuovo proved his famous contraction
criterion, establishing that any
birational contraction of smooth
surfaces contracts rational curves
of self-intersection (-1). This is when
the theory of minimal surfaces was born:
a surface is minimal if it cannot be
contracted. Jointly with his student
Enriques, Castelnuovo classified 
minimal surfaces in a series of works
which extended into 1920-ies.

Kodaira extended the minimal model theory
from projective surfaces to complex surfaces,
showing that most results of Castelnuovo and
Enriques apply also to compact complex surfaces.

The minimal model program, due to Shigefumi 
Mori, generalizes the Castelnuovo-Enriques
minimal model theory from surfaces to projective
manifolds in higher dimension. Its version
for K\"ahler threefolds was more recently
completed by H\"oring and Peternell
(\cite{_Hoering_Peternell:bir_survey_}).  However,
for non-K\"ahler complex manifolds in dimension $>2$, not much is known
about their bimeromorphic geometry.

The simplest species of non-K\"ahler manifolds
are LCK (locally conformally K\"ahler) manifolds.
These are complex manifolds which admit a K\"ahler
covering $\tilde M \arrow M$, with the deck transform group acting
by homotheties. For compact $M$, this geometry is
sharply distinct from the K\"ahler geometry. Indeed,
by Vaisman's theorem (\cite{va_tr}),
a compact LCK manifold does not admit a K\"ahler metric,
unless the deck transform acts on $\tilde M$
by isometries.

Surprizingly, the minimal model program
for LCK manifolds (at least, for a significantly
large class of LCK manifolds) is much simpler
comparing with the K\"ahler or even projective manifolds.
Let $M'$ be a compact complex
variety bimeromorphic to an LCK manifold  $M$
which belongs to two major subclasses of 
LCK manifolds, LCK manifolds with potential
or OT manifolds (Subsection \ref{_LCK_pot:Subsection_}). Then the map
$M' \dashrightarrow M$ is holomorphic.

In other words, $M$ has a unique minimal model.

The proof is complex-geometric in nature,
not using any of the state of the art 
results of birational and bimeromorphic geometry. 
However, it seems to be possible to deduce this 
result from the weak factorization theorem 
(\cite[Theorem 5-1-1]{_Matsuki:Lectures_}).

Our proof works even in a bigger generality:
let $M$ be a compact complex manifold,
and $\tilde M$ its covering. If
$\tilde M$ admits a K\"ahler form which is
exact, 
then any bimeromorphic map
$M' \dashrightarrow M$ is holomorphic
(\ref{_more_gene_Theorem_}).

\section{Preliminaries}

We present briefly the notions we need from locally conformally K\"ahler geometry. For more details examples and an up to date account of the results, please see \cite{_OV:book_}. 

\subsection{LCK manifolds}

Locally conformal K\"ahler manifolds were defined by Izu Vaisman in \cite{va_isr}. In this paper, we shall not need the orginal definition, but the following equivalent one:

\hfill

\definition (\cite[Remark 2.9]{va_tr}) A Hermitian manifold  $(M,I,\omega)$ is called {\bf Locally Conformally K\"ahler} if it admits a K\"ahler cover $(\tilde M, I, \tilde \omega)$ with deck group $\Gamma$ acting by homotheties with respect to the K\"ahler metric.
	
	The Hermitian form $\omega$ is then called {\bf an LCK form}. 
	
\hfill

\definition \label{_Character:Definition_} Let $(M,I,\omega)$ be LCK. The very definition implies the existence of a group morphism $\chi:\  \pi_1(M)\to\R^{>0}$, given by $\chi(\gamma)=\frac{\gamma^*\tilde\omega}{\tilde\omega}$ for each $\gamma\in \pi_1(M)$ viewed as a deck transformation of the K\"ahler universal cover. This group morphism is called {\bf the homothety character}.

\hfill

\definition\label{_Minimal_cover:Definition_} Let $(M,I,\omega)$ be an LCK manifold and $\chi$ the associated homothety character. The {\bf minimal K\"ahler cover} of $(M,I,\omega)$ is the K\"ahler cover associated to $\ker\chi\subset\pi_1(M)$.  

\hfill

\example \label{_Hopf_Example_}
Let $\gamma:\; \C^n \arrow C^n$ be an invertible holomorphic contraction
with apex in 0. Then the quotient $H=\frac{\C^n\setminus 0}{\langle \gamma\rangle}$
is called {\bf a Hopf manifold}.  All  Hopf manifolds are LCK. For the case
$\gamma\in\GL(n,\C)$ (when we speak about ``linear Hopf
manifolds''), the proof was given gradually, in a series
of papers: \cite{go, ov_lckpot, ov_pams}; for non-linear
contractions, the proof appeared only recently, in
\cite{ov_non_linear}.

\hfill

\example Almost all the Inoue surfaces (\cite{inoue}) are LCK (\cite{tric, bel}). The Oeljeklaus-Toma (OT) manifolds, which are higher dimensional generalizations of the Inoue surface $S^0$ are LCK. For details, see \cite[Chapter 22]{_OV:book_}.

\subsection{LCK manifolds with potential}\label{_LCK_pot:Subsection_}
	
\definition Let $(M,I,\omega)$	be an LCK manifold. It is called {\bf LCK with potential} if:
\begin{description}
	\item[(i)] The K\"ahler form $\tilde \omega$ has a smooth, 
positive K\"ahler potential: $\tilde \omega=dd^c\phi$, with $\phi:\tilde M\arrow \R^{>0}$, and
	\item[(ii)] The deck group acts by positive homotheties with respect to the potential: $\gamma^*\phi=c_\gamma\phi$, with $c_\gamma\in \R^{>0}$, for all $\gamma\in \Gamma$.
\end{description}
By abuse of terminology, we say that ``$\omega$ has potential $\phi$''. 

\hfill

\remark \label{_Automorphic_form_definition:Remark_} In general, a differential form $\eta$ on $\tilde M$ with the property that $d^*\eta=c_\gamma\eta$, $c\in\R$, $\gamma\in\Gamma$, is called {\bf automorphic}. In particular, for an LCK manifold, the K\"ahler form $\tilde\omega$ on a K\"ahler covering is always automorphic.

\hfill

\remark 
 The K\"ahler forms on the universal covers of Inoue surfaces $S^0$ and the Oeljeklaus-Toma manifolds do have positive, global potentials. However,  these are not automorphic in the  sense of \ref{_Automorphic_form_definition:Remark_} (see \cite{ot}, \cite[Chapter 22]{_OV:book_}). 

\hfill

\proposition\label{_Submanifolds_in_LCK:Proposition_} 
All smooth submanifolds of an LCK manifold with potential are LCK manifolds with potential. \endproof

\hfill

\example All Hopf manifolds (\ref{_Hopf_Example_})
are LCK with potential (\cite{ov_lckpot, ov_pams, ov_non_linear}). 
All non-K\"ahler elliptic surfaces are LCK with potential 
(\cite{bel,_ovv:surf_}). 

%

\hfill 

\ref{_Submanifolds_in_LCK:Proposition_} admits the following partial converse:

\hfill

\theorem (\cite{ov_lckpot,ov_indam})
Let $(M,I,\omega)$
be a compact LCK manifold with potential, 
$\dim_\C M\geq3$. Then $(M,I)$ admits a holomorphic embedding to a
linear Hopf manifold. \endproof

\hfill

In conclusion, if we restrict to complex dimension at
least 3, we can say that compact LCK manifolds with
potential are smooth submanifolds of linear Hopf
manifolds. If the Global Spherical Shell (GSS) conjecture (also called the Kato
conjecture) is true, the same holds for dimension 2
(\cite{ov_indam}).

\section{Bimeromorphisms of LCK manifolds}

\subsection{Normal varieties}

We shall need a result about normal varieties in the analytic category.
Recall that a complex variety $X$ is called {\bf normal}
if any locally bounded meromorphic function on an open subset 
$U\subset X$ is holomorphic 
(\cite[Definition II.7.4]{demailly}). 

\hfill

\proposition\label{_normal_bimero_Proposition_} 
Let $Z$ be a normal variety, and
$\phi:\; Z_1 \arrow Z$ a holomorphic, 
closed map such that $\phi^{-1}(z)$ is finite for all
$z$ and bijective in a general point. 
Then $\phi^{-1}$ is holomorphic.

\hfill

\proof
For proper morphisms of
algebraic varieties, this statement serves
as one of the definitions of normality: $Z$ is
normal if any finite, birational, regular map
$Z_1 \arrow Z$ is an isomorphism.

When $f$ is bijective, this statement can be found in
\cite[Prop. 14.7]{_Remmert_} or in 
\cite[Theorem 1.102]{_Greuel_Lossen_Shustin_}.

By \cite[Lemma 1.54]{_Greuel_Lossen_Shustin_}, for
any proper map $\phi:\; Z_1 \arrow Z$ such that
$\phi^{-1}(z)$ is always finite, there exist
an open neighbourhood $U_z$ for each $z\in Z$
such that $\phi^{-1}(U_z)$ is a disjoint 
union of open subsets $V_1, ..., V_n$; each of 
these would give a coordinate neighbourhood
for some $z_i\in \phi^{-1}(z)$.

Since $\phi$ is bijective in a general point,
the number of open subsets obtained by 
application of this lemma is just 1;
this implies that $\phi$ is bijective, and 
\cite[Theorem 1.102]{_Greuel_Lossen_Shustin_}
can be applied. \endproof 

%
%
%


\subsection{Fundamental group and holomorphic maps}

\proposition\label{_bimero_same_pi_1_Proposition_}
Let $\phi:\; M \dashrightarrow M_1$ be a bimeromorphic map 
of compact complex connected manifolds, and $X\subset M\times M_1$
the closure of its graph.\footnote{By definition of a meromorphic
morphism, $X$ is a complex subvariety of $M \times M_1$.} 
Then the natural projections $X \arrow M$ and
$X \arrow M_1$ induce isomorphisms of the
fundamental groups. 

\hfill

\pstep
 Denote by $\tilde X$ the
resolution of singularities of $X$.
By \cite[\S 7.8.1]{_Kollar:Shafarevich_},
the bimeromorphic holomorphic maps
$\tilde X \arrow M_1$ and $\tilde X \arrow M$
induce isomorphisms on the fundamental groups.

\hfill

{\bf Step 2:} We show that the variety $X$ is normal.
By definition, a normal variety is one where
all locally bounded meromorphic functions are holomorphic.
Note that the projections from $X$ to $M$ and
to $M_1$ are bimeromorphic; this allows us
to interpret the meromorphic functions on $X$
as meromorphic functions on $M$ and $M_1$.
Any locally bounded meromorphic function
$f$ on $X$ defines a locally bounded meromorphic function
on the manifolds $M$ and $M_1$, which are smooth and hence normal.
Therefore, $f$ is holomorphic on $M$ and on $M_1$.
This implies that $f$ is the pullback 
of a holomorphic function on $M$ and on $M_1$, 
hence it is holomorphic on $X \subset M \times M_1$.

\hfill

{\bf Step 3:} The composition of the maps
$\pi_1(\tilde X) \arrow \pi_1(X) \arrow\pi_1(M)$
(respectively  $\pi_1(\tilde X) \arrow \pi_1(X) \arrow \pi_1(M_1)$) 
is an isomorphism by Step 1. Therefore, to show that
$\pi_1(X) \arrow \pi_1(M_1)$ (respectively $\pi_1(X) \arrow \pi_1(M)$) 
is an isomorphism, it would suffice to show that the
natural map $\pi_1(\tilde X) \arrow \pi_1(X)$
is surjective.

Let $M^\circ\subset M$ be the complement of the exceptional set
of $\phi$. Clearly, the graph 
of $\phi\restrict{M^\circ}$ is dense and Zariski open in $X$.

Let $U\subset Z$ be a Zariski open subset
in a normal complex variety. Then the natural map $\pi_1(U) \arrow \pi_1(Z)$
is surjective (\cite[IX, Cor. 5.6]{_Grothendieck:SGA1_}, 
\cite[\S 1.3]{_Campana:twistors_}, or 
\cite[Lemma 3.3]{_Kollar:fund_rat_con_}).

We apply this argument to $M^\circ\subset X$ and 
$M^\circ\subset \tilde X$, and obtain
that the natural maps 
$\pi_1(M^\circ) \arrow  \pi_1(\tilde X)$ and $\pi_1(M^\circ) \arrow  \pi_1(X)$ 
are surjective. Then the natural map
$\pi_1(\tilde X) \arrow \pi_1(X)$ is also surjective.
\endproof

\subsection{Manifolds bimeromorphic to an LCK manifold}

Let $\chi:\; \pi_1(M) \arrow \R^{>0}$
be the homothety character of the LCK 
structure on $M$ (\ref{_Character:Definition_});
it is a homomorphism from the deck transform group of $\tilde M$
to $\R^{>0}$ taking a map $\Psi:\; \tilde M \arrow \tilde
M$ to the scale factor $\frac{\Psi^*\tilde \omega}{\tilde \omega}$.

For our main result, we use the following proposition.

\hfill

\proposition\label{_p_1_subvariety_Proposition_}
Let $(M, \omega)$ be a compact LCK manifold,
$(\tilde M, \tilde \omega)$ its minimal
K\"ahler cover (\ref{_Minimal_cover:Definition_}), and $Z\subset M$ a subvariety of positive
dimension. Assume that the
K\"ahler form $\tilde \omega$ is exact.
Then the image of $\pi_1(Z)$ in $\pi_1(M)$ 
contains an infinite cyclic subgroup.

\hfill

\proof 
Denote by $\tilde Z\subset \tilde M$ the 
cover of $Z$ obtained by the homotopy lifting lemma. 
If the image of $\pi_1(Z)$ in $\pi_1(M)$ 
is finite, the variety $\tilde Z$ is compact.
This is impossible, because $\tilde Z$ admits a K\"ahler
form $\tilde \omega$  which is exact, hence 
$0=\int_{\tilde Z} \tilde \omega^{\dim_\C Z}=\Vol(\tilde Z)>0$;
a contradiction. By the same argument,
the K\"ahler form $\tilde \omega$ restricted to $\tilde Z$
is not the pullback of a K\"ahler form on $Z$.
This implies that the deck transform group acts on
$(\tilde Z, \tilde \omega\restrict{\tilde Z})$
by non-trivial homotheties, implying that 
$\chi(\pi_1(Z))\subset \R^{>0}$ is non-trivial. 

Consider an element
$\gamma\in \pi_1(Z)$ such that $\chi(\gamma)\in
\R^{>0}\backslash 1$.
Then $\gamma$ is of infinite order; its image
in $\pi_1(M)$ is also of infinite order, because
$\chi$ is factorized through $\pi_1(M)$.
\endproof

\hfill

\theorem\label{_bimero_to_LCK_Theorem_}
Let $M$, $M_1$ be compact complex manifolds
and $\phi:\; M_1\dashrightarrow M$ a bimeromorphism.
Assume that $M$ is an LCK manifold, and
$(\tilde M, \tilde \omega)$ its minimal K\"ahler cover;
assume also that the  K\"ahler form $\tilde \omega$ on $\tilde M$ is
exact. Then $\phi$ is holomorphic.

\hfill

\pstep
Let $X \subset M\times M_1$ be the graph of $\phi$.
By definition, $X$ is a complex subvariety of $M\times M_1$
which projects to $M$ and $M_1$ bijectively in a general
point. We denote by $\sigma:\; X \arrow M$,
$\sigma_1:\; X  \arrow M_1$ the projection maps.
To prove \ref{_bimero_to_LCK_Theorem_},
we need to show that $\sigma_1^{-1}(z)$
is finite for all $z\in M_1$. 
Then \ref{_bimero_to_LCK_Theorem_} follows 
from \ref{_normal_bimero_Proposition_},
because $M_1$ is smooth, and therefore normal.

\hfill

{\bf Step 2:}
Assume, on the contrary, that for some $z\in M_1$,
its preimage $Z_1:=\sigma_1^{-1}(z)$ is positive-dimensional.
Since the projection of $M \times M_1$ to $M$ is
bijective on the set $M \times \{z\}$,
the set $Z_1$ projects to $M$ holomorphically
and bijectively. Let $Z\subset M$ be the image of $Z_1$ in $M$.
By Remmert's proper mapping theorem 
(\cite[\S 8.2]{demailly}), $Z$ is
a complex subvariety in $M$.

\hfill

{\bf Step 3:} 
By \ref{_p_1_subvariety_Proposition_},
the image of $\pi_1(Z)$ in $\pi_1(M)$
contains an infinite order cyclic subgroup.
Therefore, its image in $\pi_1(X)=\pi_1(M)$ also 
contains an infinite order cyclic subgroup.
This is impossible, because $\pi_1(X)= \pi_1(M)= \pi_1(M_1)$,
and the projection of $Z$ to $M_1$ is a point.
\endproof	

\hfill

Since the above proof did not use the full strength of LCK geometry
the same argument can be used to prove the following
result, which might be independently useful.

\hfill

\theorem\label{_more_gene_Theorem_}
Let $M$, $M_1$ be compact complex manifolds
and $\phi:\; M_1\dashrightarrow M$ a bimeromorphism.
Assume that $M$ admits a cover  which
admits an exact K\"ahler form. 
Then $\phi$ is holomorphic.
\endproof

\hfill

{\bf Acknowledgements:} We thank Florin Ambro and Marian
Aprodu  for very useful discussions and for
bibliographical hints. We are most of all grateful to Victor
Vuletescu for finding an error in an early version of
the paper.

\hfill

{\scriptsize

	\noindent {\sc Liviu Ornea\\
		{\sc University of Bucharest, Faculty of Mathematics and Informatics, \\14
			Academiei str., 70109 Bucharest, Romania}, and:\\
		Institute of Mathematics ``Simion Stoilow" of the Romanian
		Academy,\\
		21, Calea Grivitei Str.
		010702-Bucharest, Romania}\\
	{\tt lornea@fmi.unibuc.ro,   liviu.ornea@imar.ro}
	
	\hfill

	\noindent
	{\sc Misha Verbitsky\\
		{\sc Instituto Nacional de Matem\'atica Pura e
			Aplicada (IMPA) \\ Estrada Dona Castorina, 110\\
			Jardim Bot\^anico, CEP 22460-320\\
			Rio de Janeiro, RJ - Brasil }\\
		also:\\
		Laboratory of Algebraic Geometry, \\
		Faculty of Mathematics, National Research University 
		Higher School of Economics,
		6 Usacheva Str. Moscow, Russia}\\
	\tt verbit@verbit.ru, verbit@impa.br

}

\end{document}